\documentclass[11pt]{article}

\title{\sc\Large Localized cohomology and some applications of Popa's cocycle superrigidity theorem}
\author{Asger T\" ornquist}

\usepackage{amsmath}
\usepackage{amstext}
\usepackage{amssymb}
\usepackage{amsthm}
\usepackage{enumerate}
\usepackage{mathrsfs}

 \DeclareMathOperator{\Aut}{Aut}
 
 \DeclareMathOperator{\Inn}{Inn}

\DeclareMathOperator{\dom}{dom}

\DeclareMathOperator{\ran}{ran}

\DeclareMathOperator{\proj}{proj}

\DeclareMathOperator{\id}{Id}
\DeclareMathOperator{\Char}{Char}\DeclareMathOperator{\TFA}{TFA}
\DeclareMathOperator{\Abel}{ABEL}

\newcommand{\act}{\mathcal A}

\def\N{{\mathbb N}}
\def\Z{{\mathbb Z}}

\def\F{{\mathbb F}}
\def\P{{\mathbb P}}
\def\T{{\mathbb T}}

\def\({{\normalfont (}}
\def\){{\normalfont )}}

\def\pT{{\normalfont (T)}}

\newcounter{mypar}
\newcounter{thmcounter}

\newcommand{\mysec}[1]{\setcounter{thmcounter}{0}\addtocounter{mypar}{1}\section*{\begin{center}\normalsize{\sc \S \arabic{mypar}.
#1}\end{center}}}

\newcommand{\mythm}[2]{\addtocounter{thmcounter}{1}\subparagraph{{\sc \arabic{mypar}.\arabic{thmcounter}. #1}}{\it
#2}}

\newcommand{\thm}[2]{{\sc#1\ }{\it
#2}}

\newcommand{\remark}[1]{\subparagraph{\it #1}}

\newcommand{\defn}{\addtocounter{thmcounter}{1}\subparagraph{{\rm \arabic{mypar}.\arabic{thmcounter}.}
{\it Definition.}}}

\newcommand{\myheadpar}[1]{\addtocounter{thmcounter}{1}\subparagraph{\rm \arabic{mypar}.\arabic{thmcounter}. {\it #1}}}

\newcommand{\claim}[1]{\subparagraph{{\sc #1}}}

\begin{document}

\maketitle

\begin{abstract}
We prove that orbit equivalence of measure preserving ergodic a.e.
free actions of a countable group with the relative property (T) is
a complete analytic equivalence relation.
\end{abstract}

\mysec{Introduction}

In \cite{popa}, S. Popa introduced the notion of quotients of
Bernoulli shifts in order to obtain an infinite family of measure
preserving ergodic a.e. free orbit-inequivalent actions of a
countable group $\Gamma$ with the relative property (T) over an
infinite normal subgroup $\Lambda$. These actions are defined as
follows: Let $A$ be a countable Abelian group, and let $\hat A$ be
its dual (character) group, equipped with the normalized Haar
measure. Let $X=\hat A^\Gamma$, equipped with the product measure.
Then the (left) shift-action of $\Gamma$ on $\hat A^\Gamma$ commutes
with the action of $\hat A$, and we obtain a measure preserving a.e.
free ergodic action $\sigma^{\hat A}$ of $\Gamma$ on the quotient
$\hat A^\Gamma/\hat A$.

Popa proved in \cite{popa} that $\langle\sigma^{\hat A}: A \text{ is
torsion free countable abelian}\rangle$ is a family of
$\Gamma$-actions that are orbit equivalent precisely for isomorphic
groups, when $\Gamma$ is a group with the relative property (T) over
an infinite normal subgroup $\Lambda$. The first aim of this paper
is to prove this without any normality assumption on the subgroup
$\Lambda$:

\medskip

\thm{Theorem 1.}{Suppose $\Gamma$ is a countable discrete group with
the relative property \pT over an infinite subgroup $\Lambda$, and
$A$ and $A'$ are countably infinite Abelian groups. Then
$\sigma^{\hat A}$ and $\sigma^{\hat A'}$ are orbit equivalent iff
$A\simeq A'$.}

\medskip

The result relies on Popa's cocycle superrigidity Theorem,
\cite{popa2}. Specifically, we exploit the ``local'' untwisting
Theorem, \cite[Theorem 5.2]{popa2}, to obtain information about the
action of the subgroup $\Lambda\leq\Gamma$ in relation to the action
of the ambient group $\Gamma$.

\medskip

Theorem 1 has an interesting consequence for the complexity of orbit
equivalence for groups with the relative property (T). Namely, in \S
5 we will show that the family $\langle\sigma^{\hat A}\rangle$ is
Borel with respect to the parameter $A$. More precisely, let
$\act(\Gamma,[0,1])$ denote the natural Polish space of measure
preserving actions of $\Gamma$ on $[0,1]$, and let
$\Abel_{\aleph_0}$ be the natural Polish space of countably infinite
abelian groups. Let $\act^*_e(\Gamma,[0,1])$ be the subspace of
$\act(\Gamma,[0,1])$ consisting of ergodic a.e. free
$\Gamma$-actions. We will show that there is a Borel
$f:\Abel_{\aleph_0}\to\act^*_e(\Gamma,[0,1])$ with the property that
$f(A)$ and $f(A')$ are orbit equivalent if and only if $A$ and $A'$
are isomorphic. That is, there is a {\it Borel reduction} of the
isomorphism relation for countably infinite Abelian groups to orbit
equivalence for m.p. ergodic a.e. free actions of a countable
$\Gamma$ with the relative property (T). It is known by
\cite{friedmanstanley} that the isomorphism relation in
$\Abel_{\aleph_0}$ is complete analytic, from which we obtain:

\medskip

\thm{Theorem 2.}{Suppose $\Gamma$ is a countable discrete group with
the relative property \pT over an infinite subgroup $\Lambda$. Then
orbit equivalence, regarded as a subset of
$\act^*_e(\Gamma,[0,1])\times \act^*_e(\Gamma,[0,1])$, is a complete
analytic set. In particular, it is not Borel.}

\medskip

We also obtain in Corollary 5.6 the same result for conjugacy in
$\act^*_e(\Gamma,[0,1])$: Under the assumptions of Theorem 2,
conjugacy is analytic, but not Borel.

\bigskip

{\it Organization}: In \S 2 we introduce the notion of ``localized
cohomology'' which is the central tool used to distinguish the
actions $\sigma^{\hat A}$ up to orbit equivalence. In $\S 3$ we do a
preliminary analysis that proves that a countable group $\Gamma$
with the relative property (T) has continuum many orbit inequivalent
actions. In \S 4 we refine this analysis to prove Theorem 1 above.
Theorem 2 is proved in \S 5.

\bigskip

Research for this paper was supported in part by the Danish Natural
Science Research Council grant no. 272-06-0211.

\mysec{Localized cohomology}

Let $\Gamma$ be a countable group and $\sigma$ a probability measure
preserving (p.m.p.) $\Gamma$-action on standard Borel probability
space $(X,\mu)$. Recall that a 1-cocycle for $\sigma$ is a
measurable map $\alpha: \Gamma\times X\to \T$ such that
$$
\alpha(\gamma_1\gamma_2,x)=\alpha(\gamma_1,\sigma(\gamma_2)
(x))\alpha(\gamma_2,x)\ \  (\gamma_1,\gamma_2\in \Gamma,
\mu\text{-a.e. } x\in X).
$$
The set of all such cocycles is denoted $Z^1(\sigma)$, and forms a
Polish group under pointwise multiplication, when given the subspace
topology inherited from $L^{\infty}(X, \T)^\Gamma$. A 1-{\it
coboundary} is a cocycle $\beta\in Z^1(\sigma)$ of the form
$$
\gamma_f(g,x)=f(x)^*f(\sigma(g)(x)),
$$
where $f:X\to \T$ is a measurable map. The coboundaries form a
subgroup denoted $B^1(\sigma)$. The {\it 1st cohomology group} is
then defined as
$$
H^1(\sigma)=Z^1(\sigma)/B^1(\sigma).
$$
We now introduce the notion of a localized coboundary:

\defn Suppose $\Lambda<\Gamma$ is a subgroup. We say that $\beta\in Z^1(\sigma)$ is a
{\it $\Lambda$-local coboundary} if there is a measurable $f:X\to\T$
such that
$$
(\forall \gamma\in \Lambda)
\alpha(\gamma,x)=f(x)^*f(\sigma(\gamma)(x)),
$$
i.e. if $\beta|\Lambda$ is a 1-coboundary for $\sigma|\Lambda$. We
denote by $B^1_{\Lambda}(\sigma)$ the group of $\Lambda$-local
coboundaries. The {\it $\Lambda$-localized 1st cohomology group} is
defined as
$$
H^1_\Lambda(\sigma)=Z^1(\sigma)/B^1_\Lambda(\sigma).
$$
We can make $H^1_\Lambda(\sigma)$ into a topological group by giving
it the quotient topology.

The following relativization of a result of Schmidt's in
\cite{schmidt1}, \cite{schmidt2}, was already noted in \cite{popa}
1.6.2, though not stated in this form. See also \cite{furman}
Theorem 4.2 for a more general result along these lines.

\mythm{Proposition.}{If $\Gamma$ is a countable group with the
relative property \pT over an infinite subgroup $\Lambda<\Gamma$ and
$\sigma$ is a p.m.p. action of $\Gamma$ on a standard Borel
probability space $(X,\mu)$ such that $\sigma|\Lambda$ is ergodic,
then $B^1_\Lambda(\sigma)$ is an open subgroup of $Z^1(\sigma)$, and
$H^1_\Lambda(\sigma)$ is discrete in the quotient topology.}

\begin{proof}
It suffices to show that $B^1_\Lambda(\sigma)$ contains a
neighbourhood of the identity. Let $Q\subseteq\Gamma$ be a finite
subset and $\varepsilon>0$ such that if $\pi$ is a unitary
representation of $\Gamma$ with $(Q,\varepsilon)$-invariant vectors,
then it has a non-zero $\Lambda$-invariant vector. Suppose now that
$\alpha\in Z^1(\sigma)$ is such that
$$
\|\alpha(\gamma,x)-1\|_\infty<\varepsilon^2
$$
for all $\gamma\in Q$. Consider the unitary representation $\pi$ of
$\Gamma$ on $L^2(X)$ given by
$$
\pi(\gamma)(f)(x)=\alpha(\gamma^{-1},x)^{-1}f(\gamma^{-1}\cdot_\sigma
x).
$$
Then the constant $1$ function is $(Q,\varepsilon)$-invariant. Hence
there is a $\Lambda$-invariant non-zero $f\in L^2(X)$. Invariance
amounts to
$$
\alpha(\gamma^{-1},x)^{-1}f(\gamma^{-1}\cdot_\sigma x)=f(x),
$$
for all $\gamma\in\Lambda$, which is equivalent to
$$
f(\gamma\cdot_\sigma x)=\alpha(\gamma,x)f(x).
$$
By the ergodicity of $\sigma|\Lambda$ we have that $f(x)\neq 0$
almost everywhere. Since we also have
$$
|f(\gamma\cdot_\sigma x)|=|f(\gamma\cdot_\sigma
x)\alpha(\gamma^{-1},x)^{-1})|=|f(x)|
$$
it follows that if $\psi(x)=f(x)/|f(x)|$ then $\psi:X\to\T$ and
$$
\alpha(\gamma,x)=\psi(x)^*\psi(\gamma \cdot_{\sigma} x)
$$
for all $\gamma\in\Lambda$, thus proving that $B^1_\Lambda(\sigma)$
is open in $Z^1(\sigma)$.
\end{proof}

\myheadpar{Reduced cohomology.} Along with the localized cohomology
group we also introduce the {\it reduced} localized cohomology
group, $H^1_{\Lambda,r}(\sigma)$ as follows: Let
$B^1_{\Lambda,r}(\sigma)$ consist of all $\alpha\in Z^1(\sigma)$ of
the form
$$
\alpha(g,x)=f(g\cdot x)\beta(g,x) f(x)^*
$$
where $\beta|\Lambda\times X$ is a character (does not depend on
$x\in X$). The reduced localized cohomology group is defined as
$$
H^1_{\Lambda, r}(\sigma)=Z^1(\sigma)/B^1_{\Lambda, r}(\sigma).
$$
It is clear that if we let
$$
C_{\Lambda}(\sigma)=\{\beta\in Z^1(\sigma):
(\exists\chi\in\Char(\Lambda))\beta(g,x)=\chi(g) \text{ a.e.}\}
$$
then
$$
B^1_r(\sigma)=C_{\Lambda}(\sigma)B^1(\sigma).
$$
Further, we have:

\mythm{Lemma.}{If $\sigma|\Lambda$ is weakly mixing then
$B^1(\sigma)\cap C_\Lambda(\sigma)=\{1\}$.}
\begin{proof}
It follows that for $g\in\Lambda$ we have
$$
f(g\cdot x)=\beta(g)f(x).
$$
Hence $f$ is a $\Lambda$-eigenfunction. Since the $\Lambda$-action
is weakly mixing, we must have $f=1$.
\end{proof}

\myheadpar{Local untwisting.} The notion of local untwisting of
cocycles is, of course, the crux of Popa's construction in
\cite{popa2}. Much of the point of the present paper is that local
untwisting suffices for certain applications.

Let $\Gamma$ be a countable discrete group and $\Lambda$ a subgroup,
and suppose that $\sigma$ is a p.m.p. action of $\Gamma$ on
$(X,\mu)$. We will now consider cocycles with target group $H$,
which is assumed to be in Popa's class of {\it Polish groups of
finite type}, i.e. realizable as a closed subgroup of the unitary
group of a finite countably generated von Neumann algebra. For our
purposes the reader can assume that $H$ is either countable
discrete, or is the circle group $\T$.

Recall from \cite{popa}, \cite{popa2} that an action $\sigma$ on
$(X,\mu)$ is {\it malleable} if the flip-automorphism on $X\times X$
is in the (path) connected component of the identity in the
commutator of the product action $\sigma\times\sigma$ on $X\times
X$. We will now state a ``local'' cocycle superrigidity theorem,
which was proven by Popa in \cite[Theorem 5.2]{popa2}. It plays a
key role in the arguments in this paper.

\mythm{Theorem. ({\rm ``Local'' superrigidity, S. Popa \cite{popa2}.})}{Suppose $\Lambda$
is an infinite subgroup of $\Gamma$ such that $(\Gamma,\Lambda)$ has
property \pT. Suppose $\sigma$ is a malleable p.m.p. action of
$\Gamma$ and that $\sigma|\Lambda$ is weakly mixing. If $\alpha:
\Gamma\times X\to H$ is a measurable cocycle with target group in
Popa's class, then there is a homomorphism $\rho: \Lambda\to H$ and
$\psi: X\to H$ measurable such that
$$
(\forall g\in\Lambda) \psi(g\cdot x)\alpha(g,x)\psi(x)^{-1}=\rho(g).
$$
}

\remark{Remark.} In \cite{popa2}, Popa shows that under various
additional algebraic ``weak normality'' assumptions on the group
$\Lambda<\Gamma$, the untwisting can be continued to the whole
group, thus giving a classical type superrigidity theorem.

\mythm{Corollary. }{Under the assumptions of the previous theorem,
if $\alpha\in Z^1(\sigma)$ then $\alpha|\Lambda$ is cohomologous to
a character $\chi:\Lambda\to\T$.}

\mythm{Corollary.}{Under the assumptions of the previous theorem,
$H^1_\Lambda(\sigma)$ is isomorphic to a countable subgroup of
$\Char(\Lambda)$, and $H^1_{\Lambda,r}(\sigma)=\{1\}$.}
\begin{proof}
Clear from Proposition 2.2 and the previous Corollary and the
definition of the reduced localized cohomology group.
\end{proof}

We end this section by noting a fact about localized cohomology and
how the relative property (T) ``transfers'' when we have local
untwisting of cocycles, as in Theorem 2.6. This will play a crucial
role in our arguments:

\mythm{Proposition.}{Let $\Gamma$ be a countable discrete group and
$\Lambda\leqslant\Gamma$ a subgroup. Suppose $\Gamma$ acts by p.m.p.
transformations on $(X,\mu)$ and that $\alpha: \Gamma\times X\to H$
is a measurable cocycle, and there is a homomorphism
$\rho:\Lambda\to H$ such that $\alpha|\Lambda=\rho$. If
$(\Gamma,\Lambda)$ has property \pT  then $(H,\rho(\Lambda))$ has
property \pT.}

\begin{proof}
We use Jolissaint's characterization of relative property (T), see
\cite{jolissaint}. Let $(Q,\varepsilon)$ be a Kazhdan pair for
$(\Gamma,\Lambda)$ such that any $(Q,\varepsilon)$-invariant vector
is within $\frac 1 {10}$ of a $\Lambda$-invariant vector. Let
$Q'\subseteq H$ be a finite set such that
$$
\mu(\{x\in X: \alpha(Q,x)\subseteq Q'\})>1-\frac {\varepsilon^2} 8.
$$
We claim that $(Q',\varepsilon/\sqrt 2)$ is a Kazhdan pair for
$(H,\rho(\Lambda))$. To see this, let $\pi: H\to U(\mathscr X)$ be a
unitary representation on a Hilbert space $(\mathscr X,\|\cdot\|)$
and suppose $\xi\in\mathscr X$ is a $(Q',\varepsilon/\sqrt
2)$-invariant unit vector. Define a representation $\pi^\alpha$ of
$\Gamma$ on $L^2(X,\mathscr X)$ by
$$
\pi^\alpha(g)(f)(x)=\pi(\alpha(g^{-1},x)^{-1})(f(g^{-1}\cdot x)).
$$
Then
\begin{align*}
\pi^\alpha(g_1g_2)(f)(x)&=\pi(\alpha(g_2^{-1}g_1^{-1},x)^{-1})(f(g_2^{-1}g_1^{-1}\cdot x))\\
&=\pi(\alpha(g_1^{-1},x)^{-1}\alpha(g_2^{-1},g_1^{-1}.x)^{-1})(f(g_2^{-1}g_1^{-1}\cdot x))\\
&=\pi^\alpha(g_1)(\pi^\alpha(g_2)(f))(x).
\end{align*}
Let $f_\xi(x)=\xi$ for all $x\in X$. Then for $g\in Q$ we have
\begin{align*}
\|\pi^\alpha(g)(f_\xi)-f_\xi\|_{L^2(X,\mathscr X)}^2=&\int
\|\pi(\alpha(g^{-1},x)^{-1})(f_\xi(g^{-1}\cdot x))-\xi\|^2d\mu(x)\\
=&\int_{\{x:\alpha(g,x)\in Q'\}}
\|\pi(\alpha(g^{-1},x)^{-1})(f_\xi(g^{-1}\cdot x))-\xi\|^2d\mu(x)\\
+&\int_{\{x:\alpha(g,x)\notin
Q'\}}\|\pi(\alpha(g^{-1},x)^{-1})(f_\xi(g^{-1}\cdot x))-\xi\|^2d\mu(x)\\
\leq&\frac {\varepsilon^2} 2+4\frac {\varepsilon^2} 8=\varepsilon^2.
\end{align*}
It follows that there is a $\Lambda$-invariant unit vector $f_0\in
L^2(X,\mathscr X)$ such that $\|f_0-f_\xi\|_{L^2(X,\mathscr
X)}\leq\frac 1 {10}$. Let $V_{\mathscr X}$ be the subspace of
$L^2(X,\mathscr X)$ consisting of constant functions. Since
$\|f_0-f_\xi\|\leq \frac 1 {10}$, the projection of $f_0$ unto
$V_{\mathscr X}$ is not $0$, so let $f=\proj_{V_{\mathscr
X}}(f_0)/\|\proj_{V_{\mathscr X}}(f_0)\|$ and suppose $f=f_{\xi_0}$
for some $\xi_0\in\mathscr X$. Note that $V_{\mathscr X}$ is a
$\Lambda$ invariant subspace. Since $\pi^\alpha$ is a unitary
representation, we must then have for $h\in\Lambda$ that
$$
\pi^\alpha(h)(\proj_{V_\mathscr X} f_0)=\proj_{V_{\mathscr
X}}(\pi^\alpha(h) f_0)=\proj_{V_{\mathscr X}}(f_0).
$$
Hence $f$ is $\Lambda$-invariant, and so
$\pi(\rho(h))(\xi_0)=\xi_0$. This shows that $(H,\rho(\Lambda))$ has
property \pT.
\end{proof}

\mysec{Orbit equivalence}

We consider the following set-up: $\Gamma$ is a countably infinite
group, $\sigma:\Gamma\curvearrowright (X,\mu)$ is a p.m.p. malleable
action of $\Gamma$ and $\Lambda\leq\Gamma$ is an infinite subgroup
such that $\sigma|\Lambda$ is weakly mixing. Additionally, there is
a compact 2nd countable group $K$ acting in a measure preserving way
on $(X,\mu)$, the action of which commutes with $\sigma$. The action
of $K$ gives rise to a factor $(Y,\nu)$ consisting of
$K$-equivalence classes, and we have the factor map
$$
\theta:x\to [x]_K.
$$
The measure $\nu$ is the push-forward measure of $\mu$. Note that
$(Y,\nu)$ is standard because $K$ is assumed to be compact. $\Gamma$
acts on $(Y,\nu)$ in a p.m.p. way, and we denote this action
$\sigma^K$. (The action of $K$ will always be implicit.)

The quotients of Bernoulli shifts $\sigma^{\hat A}$ discussed in \S
1 is an example of this situation. We note the following easy fact
about $\sigma^{\hat A}$:

\mythm{Lemma.}{If $\Lambda$ is an infinite subgroup of $\Gamma$,
then $\sigma^{\hat A}|\Lambda$ is mixing.}
\begin{proof}
Let $B\subseteq \hat A^\Gamma$ be Borel and $\hat A$-invariant.
Since the Bernoulli shift $\sigma$ is mixing on all infinite
subgroups it holds for all $\varepsilon>0$ that the set of
$\gamma\in\Lambda$ such that $|\mu(\sigma(\gamma)(B)\cap
B)-\mu(B)^2|\geq\varepsilon$ is finite. Hence $\sigma^{\hat
A}|\Lambda$ is mixing.
\end{proof}

The following Lemma is certainly implicit in \cite{popa}:

\mythm{Lemma.}{Let $\Gamma$ be a countable group with the relative
property \pT over an infinite subgroup $\Lambda$ and let $A$ be a
countably infinite abelian group. Suppose
$\sigma:\Gamma\curvearrowright (X,\mu)$ is a p.m.p. action with
$\sigma|\Lambda$ weakly mixing, and that $\hat A=\Char(A)$ acts on
$(X,\mu)$ in a free, measure preserving way commuting with the
action of $\Gamma$. Let $(Y,\nu)$ be the corresponding factor,
$\theta:X\to Y$ the factor map and let $\sigma^{\hat A}$ be the
quotient action. Then $H^1_{\Lambda,r}(\sigma)=\{1\}$ implies that
$H^1_{\Lambda,r}(\sigma^{\hat A})=A$.}

\begin{proof}
For each $\alpha\in Z^1(\sigma^{\hat A})$, let $\alpha'\in
Z^1(\sigma)$ be
$$
\alpha'(g,x)=\alpha(g,\theta(x)).
$$
By assumption we can find $f:X\to\T$ and $\beta\in
C_\Lambda(\sigma)$ such that
$$
\alpha'(g,x)=f(g\cdot x)\beta(g,x)f(x)^*.
$$

\claim{Claim 1.} There is a character $\chi:\hat A\to\T$ such that
$(\forall a\in \hat A) f(a\cdot x)=\chi(a)f(x)$.

\medskip

\begin{proof}[Proof of Claim 1:] To see this, note that since $\alpha'$ is $\hat A$-invariant we
have for all $a\in \hat A$ and $g\in \Lambda$ that
$$
f(g\cdot x) \beta(g, x) f(x)^*=f(g\cdot a\cdot x) \beta(g, a\cdot x)
f(a\cdot x)^*.
$$
Using that $\beta(g,x)$ does not depend on $x$ for $g\in\Lambda$,
this gives us
$$
f(g\cdot x) f(x)^*=f(g\cdot a\cdot x)f(a\cdot x)^*,
$$
and using that the $\Gamma$ and $\hat A$ actions commute this in
turn gives us
$$
f(a\cdot g\cdot x)^*f(g\cdot x)=f(a\cdot x)^* f(x).
$$
Hence $f(a\cdot x)^* f(x)$ is $\Lambda$-invariant, and since the
$\Lambda$-action is weakly mixing this means it must be constant.
Thus
$$
f(a\cdot x)=c_a f(x)
$$
for some constant $c_a$. Let $\chi(a)=c_a$.
\end{proof}

It is easy to check now that if we define
$$
\gamma_f(g,x)=f(g\cdot x)f(x)^*
$$
then this also defines a 1-cocycle of $\sigma^{\hat A}$, since
$\gamma_f(g,x)$ is $\hat A$ invariant. Moreover, $\beta(g,x)$ is
also $\hat A$-invariant, since
$$
\beta(g,x)=\alpha'(g,x)f(g\cdot x)^*f(x).
$$
Hence $\beta(g,x)$ is a $\sigma^{\hat A}$ 1-cocycle in
$C_{\Lambda}(\sigma^{\hat A})$. Let $E\subseteq Z^1(\sigma^{\hat
A})$ denote the subgroup of all 1-cocycles $\delta$ satisfying
$$
\delta(g,\theta(x))=f(g\cdot_{\sigma} x)f(x)^*
$$
for some $\hat A$-eigenfunction $f:X\to\T$. By the above we have
$Z^1(\sigma^{\hat A})=EC_{\Lambda}(\sigma^{\hat A})$, and by Lemma
2.4 we also have $E\cap C_{\Lambda}(\sigma^{\hat A})=\{1\}$, and so
it follows that
$$
H^1_{\Lambda,r}(\sigma^{\hat A})=EC_{\Lambda}(\sigma^{\hat
A})/B^1(\sigma^{\hat A})C_{\Lambda}(\sigma^{\hat
A})=E/B^1(\sigma^{\hat A}).
$$
Since $\hat A$ is compact and acts freely on $X$ it is possible for
each character $\chi: \hat A\to\T$ to find a measurable function
$f:X\to\T$ such that $f(a\cdot x)=\chi(a)f(x)$ a.e. Hence
$$
H^1_{\Lambda,r}(\sigma^{\hat A})=E/B^1(\sigma^{\hat
A})\simeq\Char(\hat A)=A.
$$
\end{proof}

Recall that if $E$ is a measure preserving equivalence relation then
$\Inn(E)$ is the group of measure preserving transformations
$T\in\Aut(X,\mu)$ such that $xET(x)$ a.e. Then we have:

\mythm{Lemma.}{Suppose $\sigma$ and $\tau$ a.e. free p.m.p. actions
of a countable group $\Gamma$ on $(X,\mu)$ generating the same orbit
equivalence relation $E_\sigma=E_\tau=E$. Suppose
$\Lambda\leqslant\Gamma$ is a subgroup and that there is
$T\in\Inn(E)$ such that $T\sigma T^{-1}|\Lambda=\tau|\Lambda$. Then
$H^1_{\Lambda,r}(\sigma)\simeq H^1_{\Lambda,r}(\tau)$.}
\begin{proof}
We may assume that $\sigma|\Lambda=\tau|\Lambda$. Let
$\alpha:\Gamma\times X\to \Gamma$ be the cocycle defined by
$\tau(\alpha(g,x))(x)=\sigma(g)(x)$. Then $\alpha|\Lambda=\id$. For
each $\beta\in Z^1(\tau)$ define
$$
\tilde\beta(g,x)=\beta(\alpha(g,x),x).
$$
Then $\beta\mapsto\tilde\beta$ is an isomorphism $Z^1(\tau)\to
Z^1(\sigma)$, since
\begin{align*}
\tilde\beta(gg',x)&=\beta(\alpha(gg',x),x)\\
&=\beta(\alpha(g,\sigma(g')(x))\alpha(g'),x),x)\\
&=\beta(\alpha(g,\sigma(g')(x)),\tau(\alpha(g',x))(x))\beta(\alpha(g',x),x)\\
&=\beta(\alpha(g,\sigma(g')(x)),\sigma(g')(x))\beta(\alpha(g',x),x)\\
&=\tilde\beta(g,\sigma(g')(x))\tilde\beta(g',x).
\end{align*}
Moreover, for $\gamma\in\Lambda$ we have
$$
\tilde\beta(\gamma,x)=\beta(\alpha(\gamma,x),x)=\beta(\gamma,x).
$$
Hence $\beta\mapsto\tilde\beta$ maps $B^1_{\Lambda,r}(\tau)$
isomorphically onto $B^1_{\Lambda,r}(\sigma)$, and so it follows
that $H^1_{\Lambda,r}(\tau)\simeq H^1_{\Lambda,r}(\sigma)$.
\end{proof}

Before stating the next Lemma, we recall various basic notions from
\cite{tornquist}. Let $E$ be a measure preserving equivalence
relation. We will say that two actions $\sigma$ and $\tau$ of a
countable group $\Gamma$ with $E_{\sigma}, E_{\tau}\subseteq E$ such
as in the previous Lemma are {\it $E$-inner conjugate on $\Lambda$}
if there is $T\in\Inn(E)$ such that
$$
T\sigma|\Lambda T^{-1}=\tau|\Lambda.
$$
Following \cite{tornquist}, we will say that a p.m.p. action
$\sigma$ of the group $\Gamma$ is {\it ergodic on $\Lambda$} (resp.
{\it weakly mixing on $\Lambda$}), where $\Lambda\leqslant\Gamma$,
just in case $\sigma|\Lambda$ is ergodic (resp. weakly mixing) as a
$\Lambda$ action. The following was proved in \cite{tornquist},
Lemma 4.1:

\mythm{Lemma.}{Suppose $\Gamma$ has the relative property (T) over
an infinite subgroup $\Lambda\leqslant\Gamma$ and let $E$ be a
measure preserving countable equivalence relation. Then there are at
most countably many ergodic on $\Lambda$ p.m.p. $\Gamma$ actions
$E_\sigma\subseteq E$ that are not $E$-inner conjugate on $\Lambda$.
}

\bigskip

With this in hand we now have:

\mythm{Theorem.}{If $\Gamma$ is a countable group with the relative
property \pT over an infinite subgroup $\Lambda$, then $\Gamma$ has
uncountably many orbit inequivalent a.e. free p.m.p. actions on a
standard Borel probability space.}

\begin{proof}
Suppose for a contradiction that there are uncountably many
non-isomorphic countably infinite groups $\langle
A_{\xi}:\xi<\omega_1\rangle$ such that $\sigma^{\hat A_{\xi}}$ (as
defined in Lemma 3.2) are orbit equivalent for all $\xi<\omega_1$.
We can assume that all $\sigma^{\hat A_\xi}$ generate the same orbit
equivalence relation $E$. By the previous Lemma it follows that
there are $\xi_1,\xi_2<\omega_1$, $\xi_1\neq \xi_2$, such that
$\sigma^{\hat A_{\xi_1}}$ and $\sigma^{\hat A_{\xi_2}}$ are
$E$-inner conjugate on $\Lambda$. But by Lemma 3.3 it then follows
that
$$
A_{\xi_1}\simeq H^1_{\Lambda,r}(\sigma^{\hat
A_{\xi_1}})=H^1_{\Lambda,r}(\sigma^{\hat A_{\xi_2}})\simeq A_{\xi_2}
$$
contradicting that $A_{\xi_1}$ and $A_{\xi_2}$ are not isomorphic.
\end{proof}

\mysec{A finer analysis}

We now aim to refine Theorem 3.5 to show that in fact the actions
$\sigma^{\hat A}$ are orbit inequivalent for non-isomorphic $A$. We
start by noting a general lemma which is interesting in its own
right:

\mythm{Lemma.}{Suppose $\Gamma$ is a countable group with the
relative property (T) over $\Lambda\leqslant\Gamma$. Suppose
$\sigma:\Gamma\curvearrowright (X,\mu)$ is an a.e. free p.m.p.
malleable action which is weakly mixing on all infinite subgroups of
$\Lambda$. Suppose $G$ is a countable group and
$\tau:G\curvearrowright (X,\mu)$ is an a.e. free p.m.p. action which
is ergodic on all infinite subgroups and such that $E_\sigma=E_\tau$
a.e. Then there is a homomorphism $\rho:\Lambda\to G$ such that
$(G,\rho(\Lambda))$ has the relative property (T),
$H_{\rho(\Lambda)}(\tau)\simeq H_{\Lambda}(\sigma)$ and
$H_{\rho(\Lambda),r}(\tau)\simeq H_{\Lambda,r}(\sigma)=\{1\}$.}

\begin{proof}
Let $E=E_\sigma=E_\tau$. Let $\alpha:\Gamma\times X\to G$ be the
cocycle defined by
$$
\tau(\alpha(\gamma,x))(x)=\sigma(\gamma)(x).
$$
Since $\sigma$ fulfills the hypothesis of the local superrigidity
Theorem 2.6, we can find $\psi:X\to G$ and a homomorphism
$\rho:\Lambda\to G$ such that
$$
(\forall \gamma\in\Lambda) \psi(\gamma\cdot_\sigma
x)\alpha(\gamma,x)\psi(x)^{-1}=\rho(\gamma).
$$

Define $\Psi(x)=\psi(x)\cdot_\tau x$. Then $\Psi\subseteq E$ and for
all $\gamma\in\Lambda$ we have
\begin{align*}
\Psi(\gamma\cdot_\sigma x)&=\psi(\gamma\cdot_\sigma x)\cdot_\tau
(\gamma\cdot_\sigma
x)=(\psi(\gamma\cdot_\sigma x)\alpha(\gamma,x))\cdot_\tau x\\
&=(\rho(\gamma)\psi(x))\cdot_\tau x=\rho(\gamma)\cdot_\tau\Psi(x).
\end{align*}

Thus $\Psi$ conjugates the $\Lambda$ and $\rho(\Lambda)$ actions via
$\rho$, that is
\begin{equation}
(\forall\gamma\in\Lambda) \Psi(\gamma\cdot_\sigma
x)=\rho(\gamma)\cdot_\tau \Psi(x).
\end{equation}

\claim{Claim 1:} $|\ker(\rho)|<\infty$.

\begin{proof}[Proof of Claim 1:]
Suppose not. The map $\Psi$ is $\ker(\rho)$ invariant by (1) and so
since $\sigma$ is ergodic on $\ker(\rho)$ by assumption, we have
that $\Psi$ is constant on a measure 1 set. But this contradicts
that $\Psi\subseteq E$.
\end{proof}

It follows that $\rho(\Lambda)$ is infinite. Since moreover
$\Psi(X)$ is $\rho(\Lambda)$ invariant, it follows by the ergodicity
assumptions for the $G$ action that $\Psi(X)$ has full measure. Let
$\Psi'$ be a Borel right inverse of $\Psi$, i.e. $
\Psi(\Psi'(y))=y$. Then $\Psi'$ is 1-1 and $\Psi'\subseteq E$, and
so $\Psi'$ is measure preserving (see \cite{kechris3}, proposition
2.1.) Thus $\mu(\Psi'(X))=1$ and so $\Psi$ is in fact a measure
preserving transformation, with $\Psi'=\Psi^{-1}$. Note that it now
follows that $\ker(\rho)=\{1\}$ so that $\rho(\Lambda)$ is in fact
isomorphic to $\Lambda$. Moreover, $(G,\rho(\Lambda))$ has property
(T) by Proposition 2.9.

\claim{Claim 2:} $H^1_{\rho(\Lambda)}(\tau)\simeq
H^1_{\Lambda}(\sigma)$ and $H^1_{\rho(\Lambda),r}(\tau)\simeq
H^1_{\Lambda,r}(\sigma)$.

\begin{proof}[Proof of Claim 2:]
The proof is similar to Lemma 3.3. After conjugating the $G$-action
with $\Psi$, we can assume that
$$
(\forall \gamma\in\Lambda) \sigma(\gamma)(x)=\tau(\rho(\gamma))(x)
\\ \text{ (a.e.)}
$$
Note that since $\Psi$ is inner, we still have that
$E_\sigma=E_\tau$. Let $\alpha_0:\Gamma\times X\to G$ be the
corresponding cocycle defined by
$\tau(\alpha_0(\gamma,x))(x)=\sigma(\gamma)(x)$. Then for
$\gamma\in\Lambda$ we have $\alpha_0(\gamma,x)=\rho(\gamma)$. Now we
can proceed exactly as in Lemma 3.3 by defining an isomorphism
$Z^1(\tau)\to Z^1(\sigma):\beta\to\tilde\beta$ by
$$
\tilde\beta(\gamma,x)=\beta(\alpha_0(\gamma,x),x)
$$
and verify that $\beta\mapsto\tilde\beta$ maps
$B^1_{\rho(\Lambda)}(\tau)$ isomorphically onto
$B^1_{\Lambda}(\sigma)$, and $B^1_{\rho(\Lambda),r}(\tau)$
isomorphically onto $B^1_{\Lambda,r}(\sigma)$.
\end{proof}
Finally $H^1_{\Lambda,r}(\sigma)=\{1\}$ follows from Corollary 2.8.
\end{proof}

We now prove the ``quotient'' version of the previous Lemma:

\mythm{Lemma}{Suppose $\Gamma$, $\Lambda$ and
$\sigma:\Gamma\curvearrowright (X,\mu)$ are as in the previous
Lemma. Suppose $A$ is a countably infinite Abelian group, $\hat
A=\Char(A)$ its dual, that $\hat A$ acts in a free, measure
preserving way on $(X,\mu)$, and that the action of $\hat A$ and
$\sigma$ commute. Let $(Y,\nu)$ be the corresponding quotient,
$\theta: X\to Y$ the quotient map, $\sigma^{\hat A}$ the quotient
action. Then if $G$ is a countable group and
$\tau:G\curvearrowright(Y,\nu)$ is a p.m.p. a.e. free action of $G$
which is ergodic on all infinite subgroups and
$E_\tau=E_{\sigma^{\hat A}}$, then there is a subgroup $K\leqslant
G$ such that $(G,K)$ has property (T) and $H^1_{K,r}(\tau)=A$.}
\begin{proof}
Since $E_{\sigma^{\hat A}}=E_\tau$ and $\sigma^{\hat A}$ and $\tau$
are a.e. free, we have a measurable cocycle $\alpha:G\times Y\to
\Gamma$ such that $\tau(g)(x)=\sigma^{\hat A}(\alpha(g,x))(x)$. Let
$\alpha':G\times X\to \Gamma$ be the lifted cocycle defined by
$\alpha'(\gamma,x)=\alpha(\gamma,\theta(x))$. Note that $\alpha'$
determines an a.e. free p.m.p. action $\tau'$ of $G$ on $X$ by
$$
\tau'(g)(x)=\sigma(\alpha'(g,x))(x)=\sigma(\alpha(g,\theta(x))(x).
$$
Namely, by this definition
\begin{equation}
\theta(\tau'(g)(x))=\sigma^{\hat
A}(\alpha(g,\theta(x))(\theta(x)))=\tau(g)(\theta(x))
\end{equation}
for all $g\in G$. Thus we have
\begin{align*}
\tau'(g_1g_2)(x)&=\sigma(\alpha'(g_1g_2,x))(x)=\sigma(\alpha(g_1,\tau(g_2)(\theta(x)))\alpha(g_2,\theta(x)))(x)\\
&=\sigma(\alpha(g_1,\tau(g_2)(\theta(x)))\sigma(\alpha(g_2,\theta(x)))(x)\\
&=\sigma(\alpha(g_1,\theta(\tau'(g_2)(x))))\tau'(g_2)(x)\\
&=\tau'(g_1)\tau'(g_2)(x).
\end{align*}
By the previous Lemma, we now have that $G$ has property (T) over
some infinite subgroup $K\leqslant G$, and that
$H^1_{K,r}(\tau')=\{1\}$. But then by the Lemma 3.2 and (2) we have
that $H^1_{K,r}(\tau)=A$, since $\tau$ and the action of $\hat A$
commute.
\end{proof}

\thm{Theorem 1.}{Suppose $\Gamma$ is a countably infinite group with
the relative property (T) over $\Lambda\leqslant \Gamma$ and $G$ is
any countably infinite group. Let $A, A'$ be countably infinite
abelian groups and let $\sigma^{\hat A}$ and $\sigma^{\hat A'}$ be
quotients of classical Bernoulli shifts of $\Gamma$ and $G$
respectively. Then if $\sigma^{\hat A}$ and $\sigma^{\hat A'}$ are
orbit equivalent, then $A$ and $A'$ are isomorphic.}

\begin{proof}
We apply the previous Lemma to $\sigma^{\hat A}$ and $\sigma^{\hat
A'}$. Then it follows that $G$ has the relative property (T) over
some infinite subgroup $K\leqslant G$ and that $A'\simeq
H^1_{K,r}(\sigma^{\hat A'})\simeq A$.
\end{proof}

\mythm{Corollary.}{Let $\Gamma$, $A,A'$ be as in Theorem 1, and let
$\sigma^{\hat A}$ and $\sigma^{\hat A'}$ be quotients of
$\Gamma$-Bernoulli shifts. Then $\sigma^{\hat A}$ and $\sigma^{\hat
A'}$ are orbit equivalent if and only if $A$ is isomorphic to $A'$.}

\begin{proof}
By Theorem 1, it suffices to note that if $A$ is isomorphic to $A'$
then clearly $\sigma^{\hat A}$ and $\sigma^{\hat A'}$ are conjugate,
so they are in particular orbit equivalent.
\end{proof}

\mysec{Orbit equivalence is not Borel}

Let $\Gamma$ be a countable group, $(X,\mu)$ a standard Borel
probability space. We denote by $\Aut(X,\mu)$ the group of all
$\mu$-measure preserving transformations of $X$, and equip it with
the weak topology (see \cite{halmos}.) We let
$$
\act(\Gamma,X)=\{\sigma\in \Aut(X,\mu)^\Gamma: (\forall g_1,g_2\in
\Gamma) \sigma(g_1g_2)=\sigma(g_1)\sigma(g_2)\}.
$$
Since this set is closed in the product topology it is Polish, and
we naturally identify $\act(\Gamma,X)$ with the space of all measure
preserving actions of $\Gamma$ on $X$. There are various natural
subspaces, namely, the a.e. free actions which we denote by
$\act^*(\Gamma, X)$ and the ergodic a.e. free actions, denoted
$\act_e^*(\Gamma,X)$.

It is natural consider the relations of conjugacy and orbit
equivalence in $\act(\Gamma, X)$, or $\act_e^*(\Gamma,X)$. We denote
them by $\simeq$ and $\simeq_{oe}$, respectively. It is easy to see
that conjugacy is, prima facie, an analytic equivalence relation
induced by the natural conjugation action of $\Aut(X,\mu)$ on
$\act(\Gamma,X)$. It can be shown (see below) that orbit equivalence
is also an analytic equivalence relation. However, Dye's Theorem
implies that orbit equivalence has only one class in
$\act_e^*(\Z,X)$, so it is certainly not just analytic here, it is
Borel. The main goal of this section is to prove:

\medskip

\thm{Theorem 2, (v.1).}{Let $\Gamma$ be a countably infinite group
with the relative property (T). Then orbit equivalence, considered
as an equivalence relation in $\act_e^*(\Gamma,X)$, is complete analytic,
and so in particular is not Borel.}

\myheadpar{Borel reducibility.} To prove Theorem 2, we will utilize
the theory of {\it Borel reducibility} of equivalence relations that
has been developed extensively in descriptive set theory. Let $X,Y$
be Polish spaces and $E$, $F$ be equivalence relations on $X,Y$,
respectively. (We do not assume that $X$ and $Y$ have any other
structure than their Polish topology, and we do not assume anything
about $E$ and $F$ for the moment, other than they are equivalence
relations.) Then $E$ is said to be {\it Borel reducible} to $F$,
written $E\leq_B F$, if there is a Borel $f:X\to Y$ such that
$$
xEy\iff f(x)F f(y).
$$
A quick introduction to the significance of this notion is given in
the introduction of \cite{tornquist}. Here it suffices to say that
$\leq_B$ gives a degree theory for the complexity of equivalence
relations.

Let $\Abel_{\aleph_0}$ denote the space of countably infinite
Abelian groups, $\simeq_{\Abel_{\aleph_0}}$ the isomorphism relation
among such groups. $\Abel_{\aleph_0}$ can be identified with the
following Polish space:
\begin{align*}
\Abel_{\aleph_0}=\{(\cdot,e)\in\N^{\N\times\N}\times \N: &((\forall
i,j,k \in\N)
(i\cdot j)\cdot k=i\cdot (j\cdot k)\wedge\\
& ((\forall j\in\N) e\cdot j=j) \wedge\\
& ((\forall
k\in\N)(\exists l\in\N) k\cdot l=e)\wedge\\
& (\forall i,j\in\N) i\cdot j=j\cdot i\}.
\end{align*}
This is clearly a closed set in the product topology, and so it is
Polish. Note that $\simeq_{\Abel_{\aleph_0}}$ is induced by the
natural action of the infinite symmetric group $S_\infty$ on $\N$.
For notational convenience, if $G\in\Abel_{\aleph_0}$ then we will
write $\cdot_G$ for multiplication in $G$ and $e_G$ for the identity
in $G$, i.e. $G=(\cdot_G,e_G)$.

It is known by Theorem 6 of \cite{friedmanstanley} that the
isomorphism relation for Abelian $p$-groups is complete analytic.
Hence Theorem 2 version 1 will follow at once from Theorem 2 version
2 below, which itself is a consequence of Theorem 1. Note that per
the usual convention in descriptive set theory,
$\simeq_{oe}^{\act_e^*(\Gamma,X)}$ denotes the restriction of
$\simeq_{oe}$ to $\act_e^*(\Gamma,X)$.

\bigskip

\thm{Theorem 2 (v.2).}{If $\Gamma$ is a countably infinite group
with the relative property \pT and $(X,\mu)$ is a standard Borel
probability space then $\simeq_{\Abel_{\aleph_0}}$ is Borel
reducible to $\simeq_{oe}^{\act_e^*(\Gamma,X)}$.}

\remark{Remark.}{Let $\TFA_{\aleph_0}$ denotes the subset of
$\Abel_{\aleph_0}$ consisting of torsion free Abelian groups. It was
shown by Hjorth in \cite{hjorth4} that if $E$ is an equivalence
relation on a Polish space $X$ and $\simeq_{\TFA_{\aleph_0}}\leq_B
E$ then $E$ cannot be Borel. Using Hjorth's technique, Downey and
Montalban have recently shown in \cite{downeymontalban} that in fact
$\simeq_{\TFA_{\aleph_0}}$ is a complete analytic subset of
$\TFA_{\aleph_0}^2$. Theorem 2, v.2 clearly shows that
$\simeq_{\TFA_{\aleph_0}}\leq_B \simeq_{oe}^{\act_e^*(\Gamma,X)}$
and so the result of Downey and Montalban gives another reason why
$\simeq_{oe}^{\act_e^*(\Gamma,X)}$ is complete analytic.

\bigskip

The proof of Theorem 2, v.2, involves an amount of coding since the
measure preserving actions we used to prove Theorem 1 are defined on
different probability spaces. We need a few general lemmata to deal
with this. The reader should know that we rely heavily on the
results in \cite{kechris}, chapters 4.F, 12 and 17 and 28; it is
indeed the correct reference for all the descriptive set theory
needed here.

\mythm{Lemma.}{Suppose $X$ is a Polish space and let $P_c(X)$ denote
the Polish space of continuous probability measures on $X$. Then
there is a Borel map $f:P_c(X)\times X\to [0,1]$ such that for all
$\mu\in P_c(X)$ the map $f(\mu,\cdot)=f_\mu$ is a $\mu$-measure
preserving bijection from a set of full $\mu$-measure in $X$ onto a
set of full measure in $([0,1],m)$, where $m$ is Lebesgue measure.}

\begin{proof}
We may assume that $X=[0,1]$. Then we can go ahead as in the proof
of \cite[17.41]{kechris}, and define
$$
f(\mu,x)=\mu([0,x]).
$$
By \cite[17.25]{kechris}, this is Borel. Exactly as in the proof of
\cite[17.41]{kechris}, we have that $f_\mu$ is a measure preserving
bijection between sets of full measure, so $f$ is as promised.
\end{proof}

Let
$$
C=\{(G,x)\in \Abel_{\aleph_0}\times \T^\N: (\forall g_1,g_2\in\N)
x(g_1\cdot_G g_2)=x(g_1)x(g_2)\}.
$$
Then for each $G\in\Abel_{\aleph_0}$ the set $C_G$ is exactly the
set of characters on the group $\langle\N,\cdot_G,e_G\rangle$. Since
$C_G$ is compact we have by \cite[28.8]{kechris} that the map
$\Char:\Abel_{\aleph_0}\to K(\T^\N)$ where $\Char(G)=C_G$ is Borel,
where $K(\T^\N)$ denotes the compact hyperspace of $\T^\N$ as
defined in \cite[4F]{kechris}. We now have

\mythm{Lemma.}{The map $H:\Abel_{\aleph_0}\to P(\T^\N)$, which
assigns to $G\in\Abel_{\aleph_0}$ the Haar measure on $\Char(G)$, is
Borel.}

\begin{proof}
Let $(O_n)$ be a countable basis for the topology on $\T$. Let $\P$
be the set of all finite partial functions $f$ with
$\dom(f)\subseteq\N$ and $\ran(f)\subseteq\N$. For each such $f$,
let
$$
U_f=\{x\in\T^\N: (\forall i\in\dom(f)) x(i)\in O_{f(i)}\}.
$$
Then $(U_f)_{f\in\P}$ forms a countable basis for the product
topology on $\T^\N$, which is invariant under the action of the full
permutation group $S_\infty$ of $\N$ on $\T^\N$. Let
$F_f=\T^\N\setminus U_f$. Define
\begin{align*}
H=\{(G,\mu)\in\Abel_{\aleph_0}\times P(\T^\N):&\mu(\Char(G))=1\wedge\\
&(\forall f\in\P)(\forall g\in \N)\mu(F_f)=\mu(F_{g.f})\}
\end{align*}
where $g.f(i)=j\iff f(g^{-1}\cdot i)=j$. By \cite[17.29]{kechris}
$H$ is Borel, and by definition we have $(G,\mu)\in H$ precisely
when $\mu$ is the Haar measure on $\Char(G)$. By the uniqueness of
the Haar measure and \cite[14.12]{kechris} it follows that $H$
defines a Borel function $\Abel_{\aleph_0}\to P(\T^\N)$ as required.
\end{proof}

If $f:X\to Y$ is Borel, $X$, $Y$ Polish spaces, and $\mu$ a measure
on $X$, then we denote by $f[\mu]$ the push-forward measure on $Y$.
(Note that our notation differs from \cite{kechris} here, but is in
line with \cite{tornquist}):

\mythm{Lemma.}{If $f:X\times Y\to Z$ is a Borel map then there is a
Borel $f^*:X\times P(Y)\to P(Z)$ such that $f^*(x,\mu)=f_x[\mu]$,
where $f_x:Y\to Z:y\mapsto f(x,y)$.}

\begin{proof}
By \cite{kechris} 17.27 and 17.40 the map $X\times P(Y)\to P(X\times
Y):(x,y)\mapsto \delta_x\times\mu$ is Borel. So by
\cite[17.28]{kechris} we have that the map $X\times P(Y)\to
P(Z):(x,\mu)\mapsto f[\delta_x\times\mu]$ is Borel. Now note that
$f[\delta_x\times\mu]=f_x[\mu]$.
\end{proof}

\begin{proof}[Proof of Theorem 2, v.2]
Let $\Gamma$ be a fixed countably infinite group and let
$X=(\T^{\N})^\Gamma$. Consider $K(X)$, the space of compact subsets
of $X$. Note that $\Gamma$ acts on $K(X)$ since it acts on $X$ by
left-shift, and for each $G\in\Abel_{\aleph_0}$, $\Char(G)=C_G$ acts
naturally on $X$. Consider the map $f:\Abel_{\aleph_0}\times X\to
K(X)$ defined by
$$
f(G,x)=[x]_{C_G}.
$$
The map $f$ is Borel since if we fix Borel $d_n:K(X)\to X$ and
$d'_n:K(\T^\N)\to \T^\N$ as in \cite{kechris} 12.13, with
$(d_n(K))_{n\in\N}$ dense in $K$ for all $K\in K(X)$ and
$(d'_n(K'))_{n\in\N}$ dense in $K'$ for all $K'\in K(\T^\N)$, then
$$
f(G,x)=K\iff (\forall n)(\exists \chi\in C_G) \chi\cdot
x=d_n(K)\wedge (\forall m) d'_m(C_G)\cdot K=K
$$
gives an analytic definition of the graph of $f$, which suffices by
\cite{kechris} 14.12. We identify the space $C_G^\Gamma/C_G$ with
the range of $f_G=f(G,\cdot)$.

Let $f^*:\Abel_{\aleph_0}\times P(X)\to P(K(X))$ be as in Lemma 5.4.
Let $H$ be as in Lemma 5.3; then we have a map
$H^\Gamma:\Abel_{\aleph_0}\to P((\T^\N)^\Gamma)$ such that
$H^\Gamma(G)$ the product measure $H(G)^\Gamma$ and this map is
Borel by (the obvious generalization of) \cite[17.40]{kechris}. Note
that $f^*(G,H^\Gamma(G))$ is the push-forward measure on
$C_G^\Gamma/C_G$ of the measure $H^\Gamma(G)$ under the map $f_G$.
Now fix a map $f_0:P(K(X))\times K(X)\to [0,1]$ as in Lemma 5.2.
Define
$$
\theta:\Abel_{\aleph_0}\times K(X)\to [0,1]:
\theta(G,K)=f_0(f^*(G,H^\Gamma(G)),K).
$$
Then for each $G\in\Abel_{\aleph_0}$ the map
$\theta_G=\theta(G,\cdot)$ defines a measure preserving bijection
between co-null subsets of $(C_G^\Gamma/C_G,f^*(G,H^\Gamma(G))$ and
$([0,1],m)$. Define $\Theta:\Abel_{\aleph_0}\times\Gamma\times
[0,1]\to [0,1]$ by
$$
\Theta(G)(g)(x)=y\iff (\exists K\in K(X))
\theta_G(K)=x\wedge\theta_G(g\cdot K)=y.
$$
Since the measure quantifiers preserve analyticity (see
\cite{kechris} p. 233) $\Theta$ is a Borel function, and by
construction $\Theta_G$ is a measure preserving $\Gamma$-action on
$[0,1]$ which is conjugate with the action of $\Gamma$ on
$C_G^\Gamma/C_G$, for all $G\in\Abel_{\aleph_0}$. Corollary 4.3 now
guarantees that $G\mapsto\Theta_G$ is a Borel reduction of
$\simeq_{\Abel_{\aleph_0}}$ to orbit equivalence in
$\act_e^*(\Gamma,[0,1])$.
\end{proof}

In order to verify Theorem 2, v.1, it suffices to prove the
following easy lemma.

\mythm{Lemma.}{If $\Gamma$ is a countable group then $\simeq_{oe}$
is an analytic subset of
$\act(\Gamma,X,\mu)\times\act(\Gamma,X,\mu)$.}
\begin{proof}
As proved in Lemma 3 in \cite{tornquist2}, there is a Borel relation
$E\subseteq \Aut(X,\mu)\times X\times X$ such that for each
$S\in\Aut(X,\mu)$ we have that
$$
\tilde S(x)=y\iff E(S,x,y)
$$
defines a measure preserving Borel function $\tilde S$ a.e. such
that $\tilde S\in S$. Define
$$
R(\sigma,x,y)\iff (\exists g\in\Gamma) E(\sigma(g),x,y).
$$
Then
$$
(\forall^\mu x,y) xE_\sigma y\iff R(\sigma,x,y)
$$
and thus
$$
\sigma\simeq_{oe}\tau\iff (\exists T\in\Aut(X,\mu))(\forall^\mu x,y)
R(\sigma,x,y)\iff R(T\tau T^{-1},x,y),
$$
which proves that $\simeq_{oe}$ is analytic, since the measure
quantifiers preserve analyticity.
\end{proof}

\remark{Remark 1.} Clearly the proof also gives a Borel reduction of
$\simeq_{\TFA_{\aleph_0}}$ to {\it conjugacy} of measure preserving
actions. We explicitly note that the following corollary, which
should be compared with the result of a similar nature for
$\Z$-actions, due to Foreman, Rudolph and Weiss, in
\cite{forrudwei}:

\mythm{Corollary.}{If $\Gamma$ is a countably infinite group with
the relative property (T) then the conjugacy relation for ergodic,
a.e. free p.m.p. actions of $\Gamma$ on $[0,1]$ is analytic, but not
Borel.}

\remark{Remark 2.} The results of \cite{tornquist} imply that under
fairly general conditions, if a countably infinite group $\Gamma$
has the relative property (T), then both conjugacy and orbit
equivalence of p.m.p. ergodic a.e. free actions of $\Gamma$ is not
classifiable by ``countable structures'' (as defined in
\cite{hjorth3}), which in particular implies that it is not possible
to Borel reduce conjugacy and orbit equivalence in this setting to
$\Abel_{\aleph_0}$. Thus we have the following:

\mythm{Corollary.}{If $\Gamma$ has the relative property (T) over an
infinite subgroup which either contains an infinite abelian
subgroup, or is normal in $\Gamma$, then
$\simeq_{\Abel_{\aleph_0}}<_B\simeq_{oe}^{\act^*_e(\Gamma,[0,1])}$.
The same holds for the conjugacy relation in
$\act^*_e(\Gamma,[0,1])$.}

\bigskip

The normality condition in Corollary 5.6 can be replaced with the
technically weaker notion of being index stable; we refer the reader
to the last section of \cite{tornquist} for details.

\medskip

{\it Note}: Since the appearance of this paper, Ioana, Kechris and
Tsankov have shown that if $\Gamma$ is any non-amenable countable
discrete group, then orbit equivalence of its measure preserving
ergodic (indeed mixing) actions are not classifiable by countable
structures, see \cite{iokets}.

\bigskip

\begin{small}

{\sc\noindent
Kurt G\"odel Research Center, University of Vienna\\
W\"ahringer Strasse 25, 1090 Vienna, Austria\\
{\it E-mail}: {\tt asger@logic.univie.ac.at}}

\medskip

\noindent {\it Formerly:} \\
\noindent {\sc Department of Mathematics, University of Toronto\\
40 St. George Street, Room 6092, Toronto, Ontario, Canada\\
}
\end{small}

\end{document}